\documentclass[reqno]{amsart}

\usepackage{amssymb}
\usepackage{graphicx}
\usepackage{url}
\usepackage{cite}

\AtBeginDocument{{\noindent\small
This is a preprint of a paper whose final and definite form will be published
in\\ Journal of Mathematical Analysis, ISSN: 2217-3412, Volume 7, Issue 1 (2016).}
\vspace{9mm}}

\begin{document}

\title[Optimal control strategies for the spread of Ebola]{Optimal
control strategies for the spread\\ of Ebola in West Africa}

\author[A. Rachah, D. F. M. Torres]{Amira Rachah, Delfim F. M. Torres}

\address{Amira Rachah \newline
\indent Math\'{e}matiques pour l'Industrie et la Physique\newline
\indent Institut de Math\'{e}matiques de Toulouse, Universit\'e Paul Sabatier\newline
\indent F-31062 Toulouse, Cedex 9, France}
\email{arachah@math.univ-toulouse.fr}

\address{Delfim F. M. Torres \newline
\indent Center for Research and Development in Mathematics and Applications (CIDMA)\newline
\indent Department of Mathematics, University of Aveiro, 3810-193 Aveiro, Portugal}
\email{delfim@ua.pt}

% ------------------------------------

\thanks{Submitted June 9, 2015. Revised and Accepted Dec 10, 2015.}

% ------------------------------------



% ------------------------------------

\begin{abstract}
The spread of Ebola virus in 2014 is unprecedented.
The epidemic is still affecting West Africa,
exacerbated by extraordinary socioeconomic disadvantages
and health system inadequacies. With the aim of understanding,
predicting, and control the propagation
of the virus in the populations of affected countries, it is crucial
to model the dynamics of the virus and study several strategies to control it.
In this paper, we present a very simple mathematical model that describes
quite well the spread of Ebola. Then, we discuss several strategies
for the control of the propagation of this lethal virus into populations,
in order to predict the impact of vaccine programmes, treatment, and
the impact of educational campaigns.
\end{abstract}

% ------------------------------------

\maketitle

\numberwithin{equation}{section}
\newtheorem{theorem}{Theorem}[section]
\newtheorem{lemma}[theorem]{Lemma}
\newtheorem{proposition}[theorem]{Proposition}
\newtheorem{corollary}[theorem]{Corollary}
\newtheorem*{remark}{Remark}

% ------------------------------------

\section{Introduction}

The first Ebola outbreak took place in 1976 in Congo, close to Ebola river,
from where the disease takes its name. Since August 2014, it is affecting
several countries in west of Africa, mainly Guinea, Sierra Leone, and Liberia
\cite{althaus,barraya,joseph,MyID:321,who}. Nowadays, Ebola is one of the most
deadliest pathogens for humans, due to the extremely rapid increase of the disease
and the high mortality rate. Several authors consider the virus a health
and humanitarian catastrophe of historic scope \cite{barraya,joseph}.

In a first stage, Ebola virus is characterised by the sudden onset of fever,
intense weakness, headache, fatigue, and muscle pain. This is followed
by vomiting blood and passive behaviour. The second stage is characterised
by diarrhea, rash, symptoms of impaired kidney and liver function, and both
internal and external bleeding (bleeding from nose, mouth, eyes and anus)
\cite{anon2,okwar,anon1}. In a third/final stage, the symptoms are summarized
in a loss of consciousness, seizures, massive internal bleeding, followed
by death \cite{borio,dowel,legrand,peter}. Ebola spreads through human-to-human
transmission via close and direct physical contact (through broken skin
or mucous membranes) with infected bodily fluids. The most infectious fluids
are blood, feces, and vomit secretions. However, all body fluids have the capacity
to transmit the Ebola. The virus is also transmitted indirectly via exposure
to objects or an environment contaminated with infected secretions. For all these
reasons, health care workers practice strict infection prevention
and control precautions when dealing with Ebola cases.

In epidemiology, mathematical models are a key tool for the understanding
of the dynamics of infectious diseases and the impact of vaccination programmes.
In fact, mathematics has an important role in the control of propagation of virus,
allowing policymakers to predict the impact of particular vaccine programmes
or to derive more efficient strategies based on mathematical insights---see
\cite{vacc_opt4,MyID:314,vacc_opt3,vacc_opt5,MyID:306,MyID:318} and references
therein. In particular, optimal control theory has become in last years
a powerful mathematical tool that can assess the intervention of public
health authorities. Indeed, the inclusion, in an epidemic model,
of some practical control strategies, like vaccines, social distancing,
or quarantine, provides a rational basis for policies, designed to
control the spread of a virus \cite{vacc_opt2,vacc_opt1}. In this spirit,
we focus our work on the investigation of effective strategies to control
the spread of the Ebola virus by setting optimal control problems subject
to a SIR epidemic model. Two practical control strategies are here considered:
vaccines or educational campaigns and treatment. The SIR model, on which our
optimal control studies are based, divides the population into three groups:
the Susceptible (S), the Infected (I), and the Recovered (R)
\cite{Gerard,chow,zeng}. This model is simulated by Rachah and Torres in their
study of the spread of the recent outbreak of Ebola virus \cite{MyID:321}.

The paper is organized as follows. In Section~\ref{sec:2}, we present a
mathematical model that describes the dynamics of the propagation of the Ebola
virus into a population. After the mathematical modelling, we use the obtained
model to discuss it in Section~\ref{sec:3} under several control strategies
for the propagation of the virus. In these strategies, we use
parameters estimated from recent statistical data based
on the WHO (World Health Organization) report of the 2014 Ebola outbreak \cite{who}.
First, we study the case of control through a vaccination strategy
(Section~\ref{subsec:3.2}); secondly, we consider educational campaigns
as a control, coupled with a treatment strategy (Section~\ref{subsec:3.3}).
Our results improve the recent investigations of \cite{MyID:321},
briefly summarized in Section~\ref{subsec:3.1}. We end with
Section~\ref{subsec:3.5} of discussion of results,
and Section~\ref{sec:4} of conclusions and future work.

% ------------------------------------

\section{Model formulation}
\label{sec:2}

\begin{figure}
\centering
\includegraphics[width=7cm]{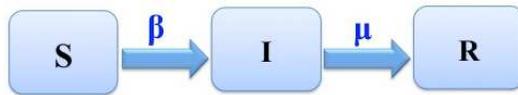}
\caption{The compartment diagram of the Susceptible--Infected--Recovered
model \eqref{eq1:SIR}.}
\label{SIR_fig1}
\end{figure}

In this section we present, mathematically, the
dynamics of the population infected by the Ebola virus.
The dynamics is described by a system of ordinary differential equations.
This system is based on the common SIR (Susceptible--Infectious--Recovered)
epidemic model, where the population is divided into three disjoint compartments:
susceptible individuals $S(t)$, i.e., people who can catch the virus;
infectious individuals $I(t)$, i.e., people who are infected by the virus
and can transmit it; and recovered individuals $R(t)$, i.e., people who have
recovered from the virus. The total population, assumed constant during the
period of time under study, is given by  $N = S(t) + I(t) + R(t)$. Ebola
spreads through human-to-human transmission via contact. The transition between
different states is described by the following two parameters: the infection
rate $\beta$ and the recovered rate $\mu$. The dynamics of the model is governed
by the following system of differential equations:
\begin{equation}
\label{eq1:SIR}
\begin{cases}
\dfrac{dS(t)}{dt} = -\beta S(t)I(t), \\[0.30cm]
\dfrac{dI(t)}{dt} =  \beta S(t)I(t)- \mu I(t),\\[0.30cm]
\dfrac{dR(t)}{dt} =  \mu I(t).
\end{cases}
\end{equation}
The first equation of the system describes the population of the susceptible
group, which reduces as the infected come into contact with them with a rate
of infection $\beta$. This means that the change in the population of susceptible
is equal to the negative product of $\beta$ with $S(t)$ and $I(t)$.
The second equation describes the infectious group over time, knowing
that the population of this group changes in two ways:
(i) people leave the susceptible group and join the infected group,
thus adding to the total population of infected a term $\beta S(t)I(t)$;
(ii) people leave the infected group and join the recovered group,
reducing the infected population by $-\mu I(t)$. The third equation
describes the recovered population, which is based on the individuals
recovered from the virus at a rate $\mu$. This means that the recovered group
is increased by $\mu$ multiplied by $I(t)$. The dynamical system \eqref{eq1:SIR}
can be represented graphically as in Figure~\ref{SIR_fig1}.

% ------------------------------------

\section{Optimal control strategies}
\label{sec:3}

Optimal control techniques are of great importance in developing optimal
strategies to control various kind of diseases and epidemics
\cite{MR2719552,delf,MR3266821,MR3101449}. For instance, to address the
challenges of obtaining an optimal vaccination strategy, one can use optimal
control theory \cite{vacc_opt4,vacc_opt3}.
In this section, we formulate three optimal control problems subject to the
SIR model \eqref{eq1:SIR}, with the aim to derive optimal
strategies. For each strategy, we study a specific objective in order
to minimize not only the number of infected individuals and the systemic costs
of vaccination or treatment, but also to include educational campaigns
within the control program. The integration of educational
campaigns has a great importance in
countries that don't have the capacity to defend themselves against the virus.
We compare the result of each strategy with the simulation results
studied by Rachah and Torres in \cite{MyID:321}, the so called Strategy~1,
which is briefly summarized in Section~\ref{subsec:3.1}. Strategies~2
and~$3$ are an improvement of Strategy~1, and are given in Sections~\ref{subsec:3.2}
and \ref{subsec:3.3}, respectively. More precisely, Strategy~3 consists in the
study of an educational campaign (about the virus) coupled with
a treatment program. The comparison between the different strategies
and the simulation results is given in Section~\ref{subsec:3.5}.

For the numerical solutions of the optimal control problems, we have used the
\textsf{ACADO} solver \cite{acado}, which is based on a multiple shooting method,
including automatic differentiation and based ultimately on the semidirect
multiple shooting algorithm of Bock and Pitt \cite{acado2}. The \textsf{ACADO}
solver comes as a self-contained public domain software environment,
written in \textsf{C++}, for automatic control and dynamic optimization.

% ------------------------------------

\subsection{Strategy 1}
\label{subsec:3.1}

In this subsection, we recall briefly the strategy studied by Rachah and Torres
in \cite{MyID:321}. We will improve it later by studying other strategies
in order to better control the propagation of the spread of Ebola into
populations. In their study, Rachah and Torres \cite{MyID:321} studied
the SIR model with control, which is given by the following system
of nonlinear differential equations:
\begin{equation}
\label{SIR_control}
\begin{cases}
\dfrac{dS(t)}{dt} = -\beta S(t)I(t) - u(t) S(t),\\[0.3cm]
\dfrac{dI(t)}{dt} = \beta S(t)I(t) - \mu I(t),\\[0.3cm]
\dfrac{dR(t)}{dt} = \mu I(t) + u(t) S(t).
\end{cases}
\end{equation}
The goal of the strategy is to reduce the infected individuals and the cost
of vaccination. Precisely, the optimal control problem consists of minimizing
the objective functional
\begin{equation}
\label{cost_func_strat1}
J(u) = \int_{0}^{t_{end}} \left[I(t) + \dfrac{\nu}{2}u^2(t)\right] dt
\end{equation}
subject to the model described by \eqref{SIR_control},
where $u(t)$ is the control variable, which represents the vaccination rate
at time $t$, and the parameters $\nu$ and $t_{end}$ denote, respectively,
the weight on cost and the duration of the vaccination program.
In our simulations we take $\beta=0.2$, $\mu=0.1$, $\nu= 0.5$,
and $t_{end}=100$ days.

% -----------------------------------

\subsection{Strategy 2}
\label{subsec:3.2}

Our goal in this strategy is to reduce the number of susceptible
and infected individuals and simultaneously increase the number
of recovered individuals. Precisely, our optimal control problem
consists of minimizing the objective functional
\begin{equation}
\label{cost_func_strat2}
J(u) = \int_{0}^{t_{end}}
\left[ A_1S(t) + A_2 I(t)
- A_3 R(t) + \dfrac{\tau}{2}u^2(t) \right] dt
\end{equation}
subject to the model described by \eqref{SIR_control}.
The two first terms in the functional objective \eqref{cost_func_strat2}
represent benefit of $S(t)$ and $I(t)$ populations
that we wish to reduce; $A_1$ and $A_2$ are
positive constants to keep a balance in the size of $S(t)$ and $I(t)$,
respectively. In our simulations we used $A_1= 0.1$ and  $A_2=0.5$.
The third term represents the recovered individuals, which we wish to increase
through vaccination ($A_3$ is a positive constant,
chosen in the simulations as $A_3=0.002$).
In the quadratic term of \eqref{cost_func_strat2}, $\tau$ is a positive
weight parameter associated with the control $u(t)$, and the square
of the control variable reflects the severity of the side effects
of the vaccination. In our numerical results we took $\tau=1$.
One has $u \in \mathcal{U}_{ad}$, where
$$
\mathcal{U}_{ad}=\left\{u : u \,  \text{is measurable}, 0
\leq u(t) \leq u_{max}<\infty, \, t\in [0,t_{end}] \right\}
$$
is the admissible control set, with $u_{max}=0.9$.
Note that in this strategy, the control variable $u(t)$ is the percentage
of susceptible individuals being vaccinated per unit of time.
Note also that in the equation $dS(t)/dt$, we have the term $-u(t)S(t)$.
Then, by minimizing the number of susceptible, the number of recovered
(described by equation $dR(t)/dt$, which depends on $u(t)S(t)$) is maximized
(through vaccination). In the numerical simulations of model \eqref{eq1:SIR}
and the strategy of the  optimal control problem, Rachah and Torres \cite{MyID:321}
used the WHO data of the 2014 Ebola outbreak occurred in Guinea, in Sierra Leone,
and in Liberia \cite{althaus,kaurov}. In order to compare our improvement
of the optimal control study with the previous results of \cite{MyID:321},
we use here the same parameters, that is, the same rate of infection $\beta=0.2$,
the same recovered rate $\mu=0.1$, and the same initial values
$\left(S(0),I(0),R(0)\right) =\left(0.95, 0.05, 0\right)$
for the initial number of susceptible, infected, and recovered populations
(at the beginning, 95\% of population is susceptible
and 5\% is infected with Ebola).

\begin{figure}
\centering
\includegraphics[width=10cm]{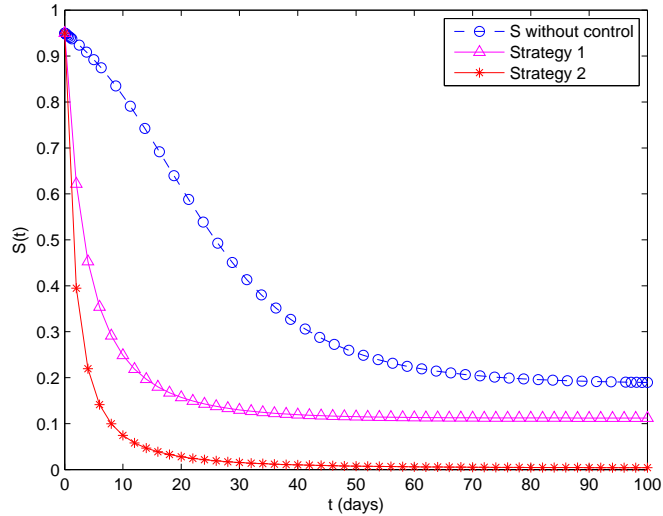}
\caption{Comparison between the curves of susceptible individuals $S(t)$
in case of Strategy~$1$ and Strategy~$2$ \emph{versus}
without control. \label{stratg12_S}}
\end{figure}

\begin{figure}
\centering
\includegraphics[width=10cm]{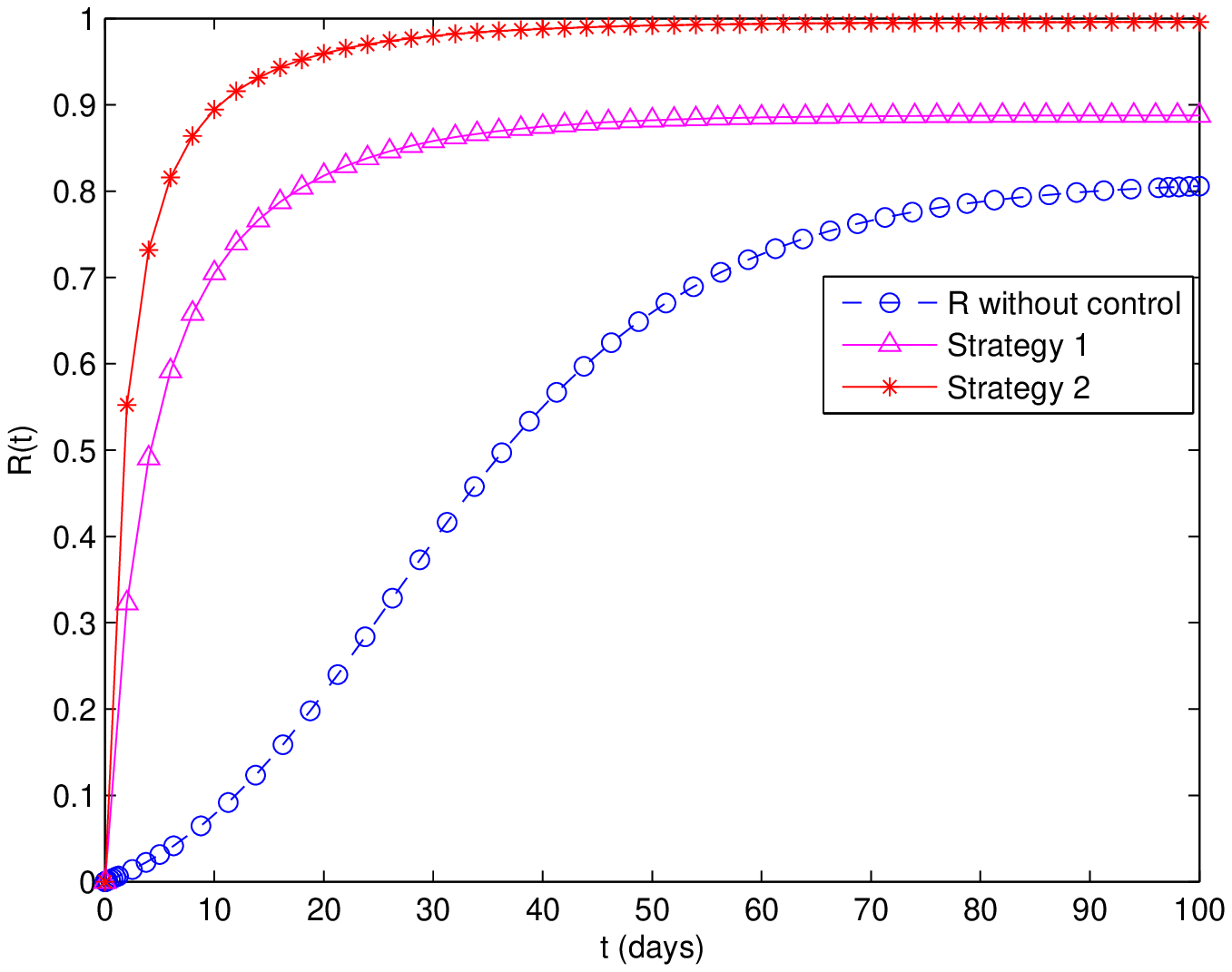}
\caption{Comparison between the curves of recovered individuals $R(t)$
in case of Strategy~$1$ and Strategy~$2$ \emph{versus}
without control. \label{stratg12_R}}
\end{figure}

\begin{figure}
\centering
\includegraphics[width=10cm]{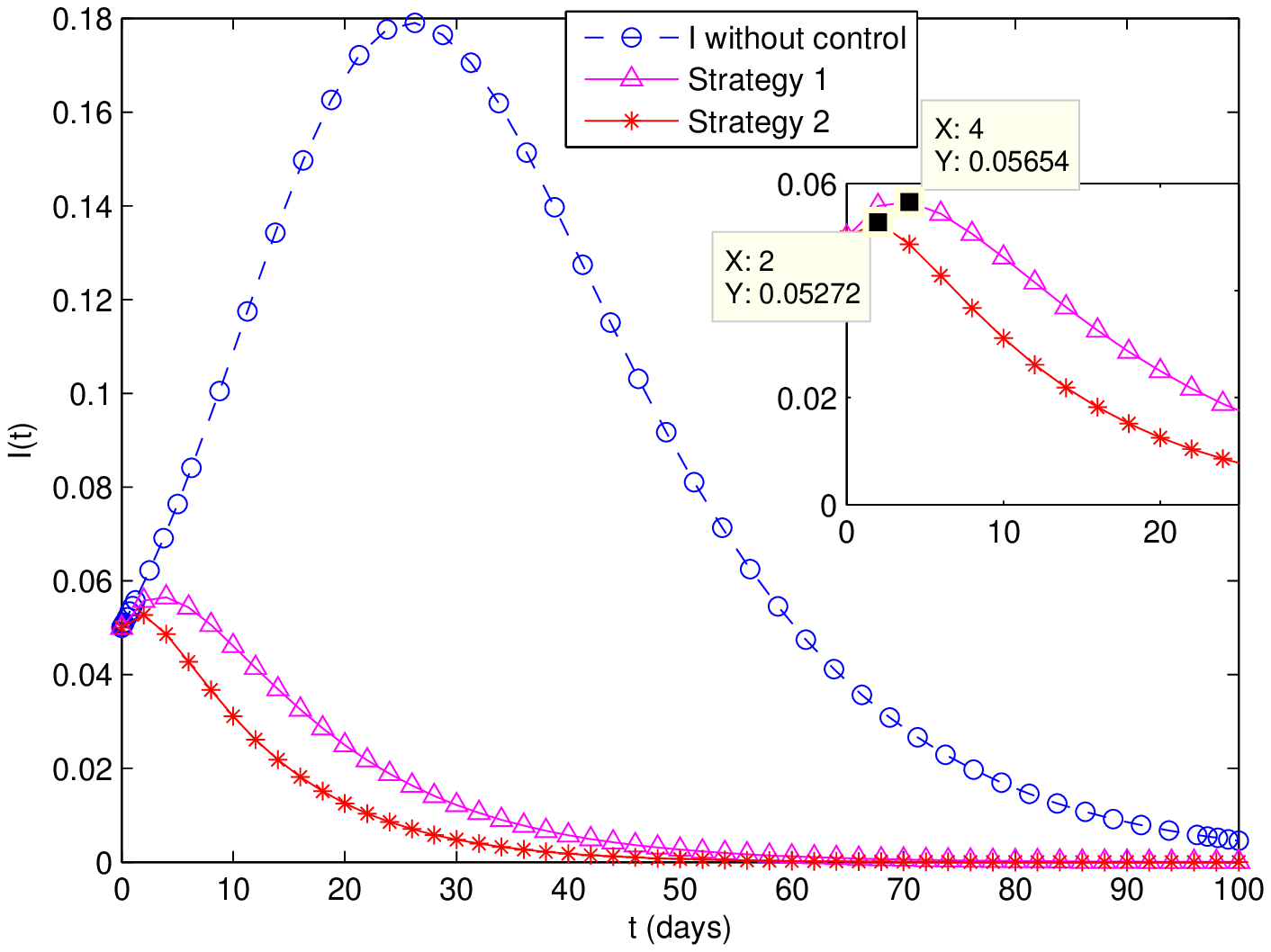}
\caption{Comparison between the curves of infected individuals $I(t)$
in case of Strategy~$1$ and Strategy~$2$ \emph{versus}
without control. \label{stratg12_I}}
\end{figure}

\begin{figure}
\centering
\includegraphics[width=10cm]{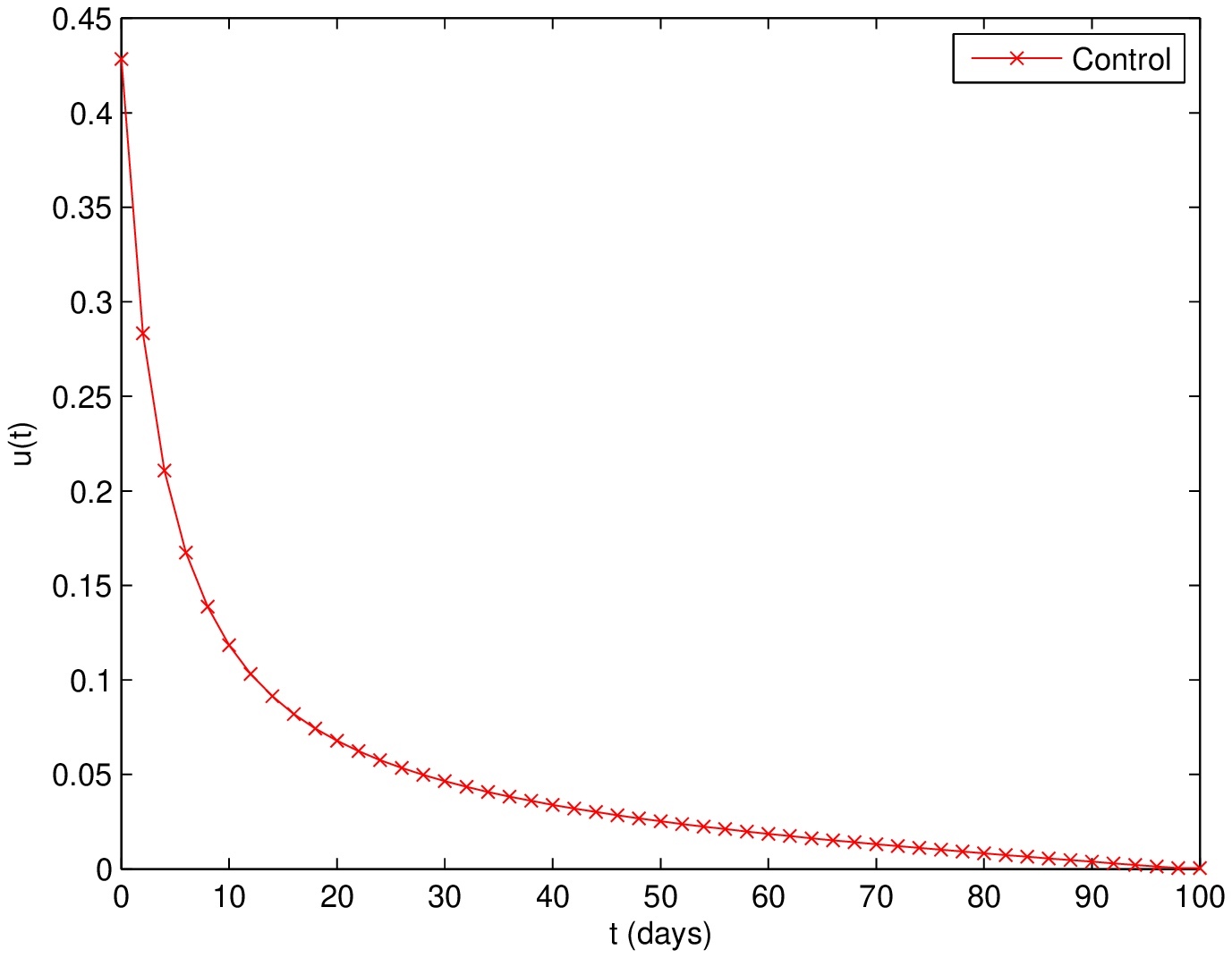}
\caption{The optimal control $u$ for Strategy~$2$
with $A_1= 0.1$, $A_2=0.5$, $A_3=0.002$, $\tau=1$
and $t_{end}=100$ days.
\label{stratg2_u}}
\end{figure}

Figures~\ref{stratg12_S}, \ref{stratg12_R} and \ref{stratg12_I} show,
respectively, the significant difference in the number of susceptible, recovered,
and infected individuals with Strategy~$1$, Strategy~$2$, and without control.
In Figure~\ref{stratg12_S}, we see that the number of susceptible $S$,
in case of optimal control under Strategy~$2$, decreases more rapidly during
the vaccination campaign. It reaches 4\% at the end of the campaign,
in contrast with the 11\% at the end of the campaign with Strategy~$1$,
and against 19\% in the absence of optimal control.
Figure~\ref{stratg12_R} shows that the number of recovered individuals
increases rapidly. The number $R(t_{end})$ of recovered at the end of the optimal
control vaccination period of Strategy~2 is 99.5\%, instead of 88.7\% in case
of Strategy~$1$, and against 80.5\% without control.
In Figure~\ref{stratg12_I}, the time-dependent curve of infected individuals
shows that the peak of the curve of infected individuals is less important
in case of control with Strategy~2. In fact, the maximum value on the infected
curve $I$ under optimal control is 5.2\% in case of Strategy~2, instead of 5.6\%
in Strategy~$1$, against 17.9\% without any control (see Figure~\ref{stratg12_I}).
The other important effect of Strategy~$2$, which we can see in the same curve,
is the period of infection, which is the less important. The value of the period
of infection is $48$ days in case of Strategy~$2$, instead of $64$ days
in case of Strategy~$1$, and against $100$ days without vaccination. This shows
the efficiency of vaccination control with Strategy~$2$ in controlling Ebola.
Indeed, the fact that in Figure~\ref{stratg12_I} the period of infection of Strategy~2
is shorter than the period of infection of Strategy~1 is important.
Figure~\ref{stratg2_u} gives a representation
of the optimal control $u(t)$ for Strategy~$2$.

% -------------------------------
% -------------------------------
% -------------------------------
% -------------------------------

% -------------------------------

\subsection{Strategy 3}
\label{subsec:3.3}

\begin{figure}
\centering
\includegraphics[scale=0.6]{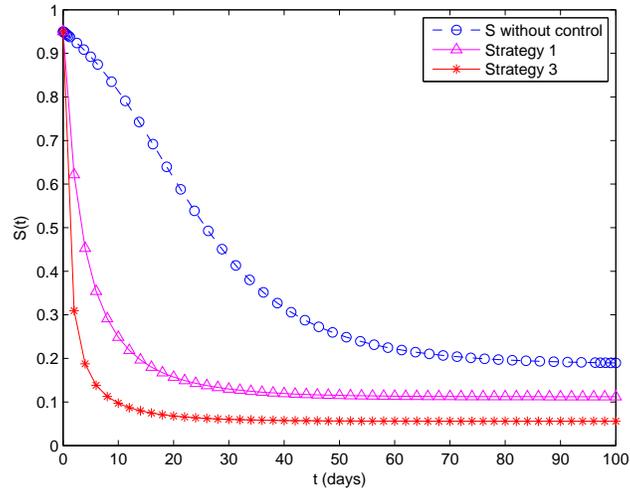}
\caption{Comparison between the curves of susceptible individuals $S(t)$
in case of Strategy~$1$ and Strategy~$3$
\emph{versus} without control. \label{stratg13_S}}
\end{figure}

\begin{figure}
\centering
\includegraphics[scale=0.6]{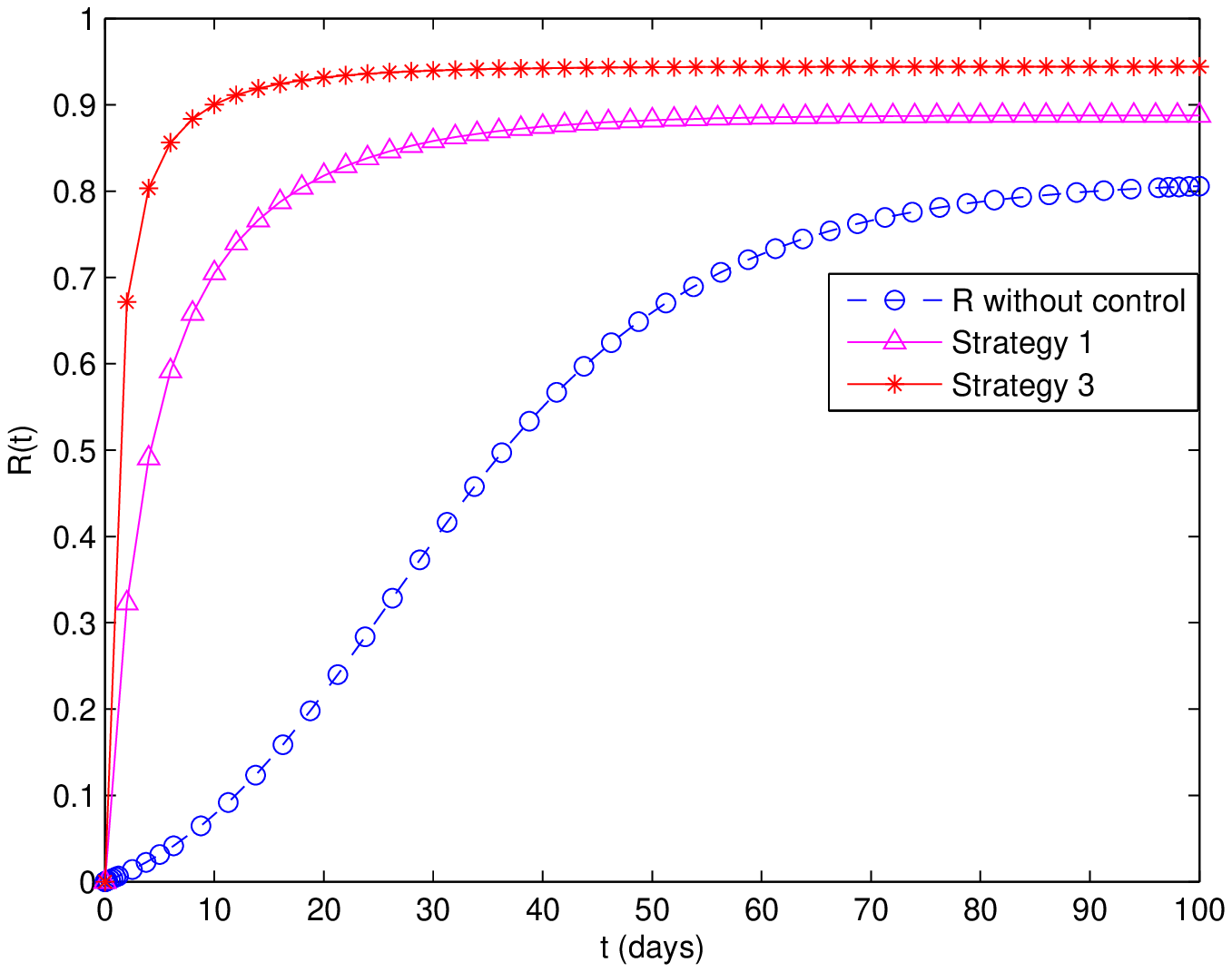}
\caption{Comparison between the curves of recovered individuals $R(t)$
in case of Strategy~$1$ and Strategy~$3$
\emph{versus} without control. \label{stratg13_R}}
\end{figure}

\begin{figure}
\centering
\includegraphics[scale=0.6]{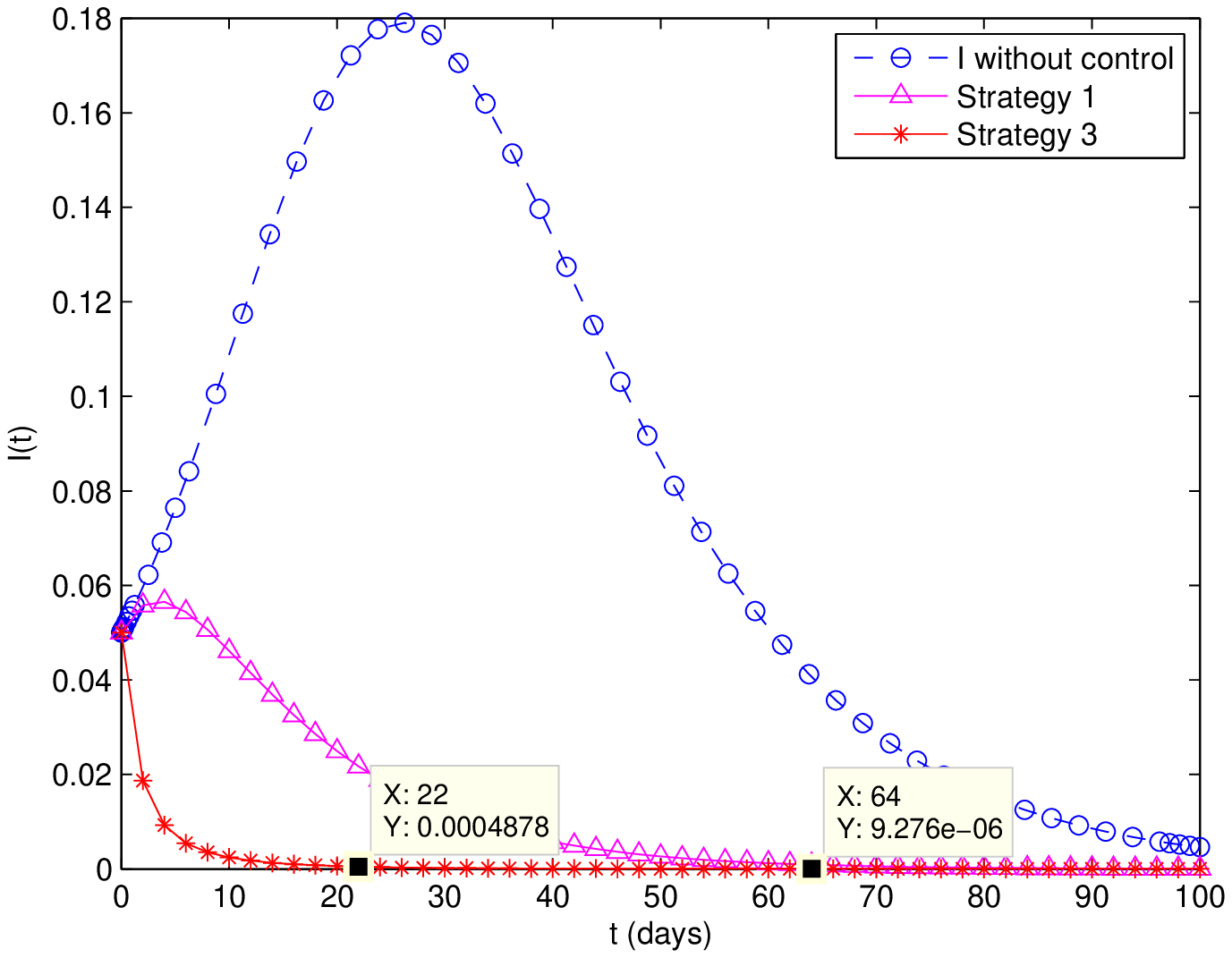}
\caption{Comparison between the curves of infected individuals $I(t)$
in case of Strategy~$1$ and Strategy~$3$
\emph{versus} without control. \label{stratg13_I}}
\end{figure}

\begin{figure}
\centering
\includegraphics[scale=0.6]{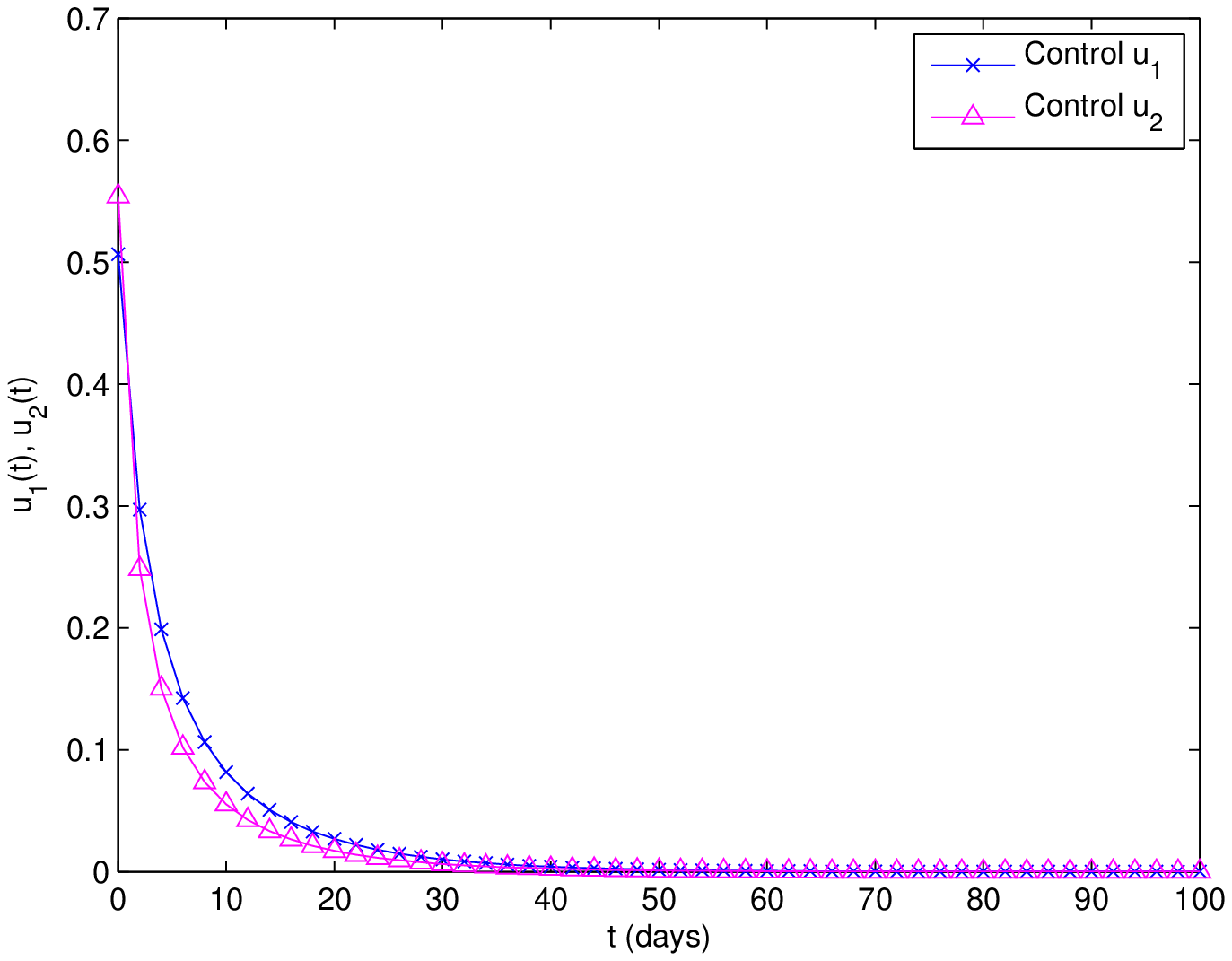}
\caption{The optimal control variables $u_1$ and $u_2$ for Strategy~$3$
with $\kappa=1$, $B_1=0.2$, $B_2=0.04$ and $t_{end}=100$ days.
\label{stratg3_u}}
\end{figure}

In this strategy we use the fact that individuals can acquire immunity against
the virus either through educational campaigns
or recovery after treatment for the virus. Our main idea is to study the effect
of educational campaigns with a treatment in practical Ebola situations. The case
of the French nurse cured of Ebola is a proof of the possibility of educational
campaigns and treatment \cite{valler}. An educational campaign, in case
of spread of Ebola, has great importance. In fact, Ebola virus spreads
through human-to-human transmission, not only by close and direct physical
contact with infected bodily fluids, but also via exposure to objects
or contaminated environment. The most infectious fluids are blood,
feces, and vomit secretions. However, all body fluids have the capacity
to transmit the virus.
Here, we intend to control the propagation of the Ebola virus by using
two control variables in the SIR model, as follows:
\begin{equation}
\label{SIR_control_stratg3}
\begin{cases}
\dfrac{dS(t)}{dt} = -\beta S(t)I(t) - u_2(t) S(t),\\[0.3cm]
\dfrac{dI(t)}{dt} = \beta S(t)I(t) - \mu I(t) -u_1(t)I(t),\\[0.3cm]
\dfrac{dR(t)}{dt} = \mu I(t) + u_1(t) I(t) + u_2(t) S(t),
\end{cases}
\end{equation}
where $u_1(t)$ is the fraction of infective that is treated per unit of time,
and $u_2(t)$  is the fraction of susceptible individuals that is subject
to an educational campaign per unit of time. Note that in Strategy~3
we study the scenario in which individuals can acquire immunity against
the virus either through educational campaigns
or recovery after treatment of the virus. The control $u_1$ (treatment)
that appears in the $dI(t)/dt$ equation is used as a treatment
applied after infection. The French nurse cured
(after infection), thanks to the treatment occurred in France, supports
this idea. Our aim is to study the application of the same treatment
(applied to the French nurse) in the countries suffering from Ebola.
Then the control $u_1(t)$ is the fraction of infected that is treated per unit of time.
Comparing with the others strategies, we would like to prove that we can apply
educational campaigns and  medical treatment (similarly to the case of the French nurse)
in the countries suffering from Ebola, which do not have the capacity
to install a vaccination program to all population. Our goal here is to minimize
simultaneously the total number of individuals that are infected, the cost
of treatment, and the cost of educational campaigns to the population.
The objective functional is now
\begin{equation}
\label{cost_func_strat3}
J(u) = \int_{0}^{t_{end}} \left[ \kappa I(t)
+ B_1\dfrac{u_1^2(t)}{2}+ B_2\dfrac{u_2^2(t)}{2} \right] dt
\end{equation}
subject to system \eqref{SIR_control_stratg3}, where
$u = (u_1, u_2)$, with $u_1$ representing treatment
and $u_2$ educational campaigns. The Lebesgue
measurable control set is defined as
$$
\mathcal{U}_{ad}:=\left\{u = (u_1,u_2) : u \,
\text{ is measurable}, \, 0 \leq u_1(t), \, u_2(t) \leq u_{max},
\, t\in [0,t_{end}] \right\},
$$
where $u_{max}=0.9$ and $\kappa$, $B_1$ and $B_2$ are weight parameters,
which in the numerical simulations we took as $\kappa=1$,
$B_1=0.2$ and $B_2=0.04$. Here, we choose quadratic terms
with respect to the controls in order to describe
the nonlinear  behaviour of the cost of implementing the educational campaigns
and the treatments. The first term in the objective functional
\eqref{cost_func_strat3}, $\kappa I$, stands for the total number
of individuals that are infected; the term $B_1 u_1^2 / 2$
represents the cost of treatment; while the term $B_2 u_2^2 / 2$
represents the cost associated with the educational campaigns.
Figure~\ref{stratg13_S} shows the time-dependent curve of susceptible
individuals, $S(t)$, which decreases more rapidly in case of the optimal control
with Strategy~$3$. It reaches 5\% at the end of the campaign, instead of 11\%
at the end of the campaign in case of Strategy~1, and against 19\%
in the absence of optimal control.
Figure~\ref{stratg13_R} shows that the number of recovered individuals
of Strategy~$3$ increases rapidly until 94.4\%, instead of 88.7\%
in case of Strategy~$1$, and against 80.5\% without control.
In Figure~\ref{stratg13_I} we see the time-dependent curve of infected
individuals $I(t)$, which decreases mostly. The other important effect
of Strategy~$3$, which we can see in the same curve, is the period
of infection, which is the less important. The value of the period of infection
is $22$ days in case of Strategy~$3$, instead of $64$ days in case of Strategy~$1$,
and against $100$ days without control. This shows the efficiency of the effect
of educational campaigns in controlling Ebola virus with the treatment
control program described in Strategy~$3$. Figure~\ref{stratg3_u} gives
a representation of the optimal control variables $u_1(t)$ and $u_2(t)$
for Strategy~$3$.
% -------------------------------
% -------------------------------
% -------------------------------
% -------------------------------

% -------------------------------

\section{Discussion}
\label{subsec:3.5}

\begin{figure}
\centering
\includegraphics[scale=0.6]{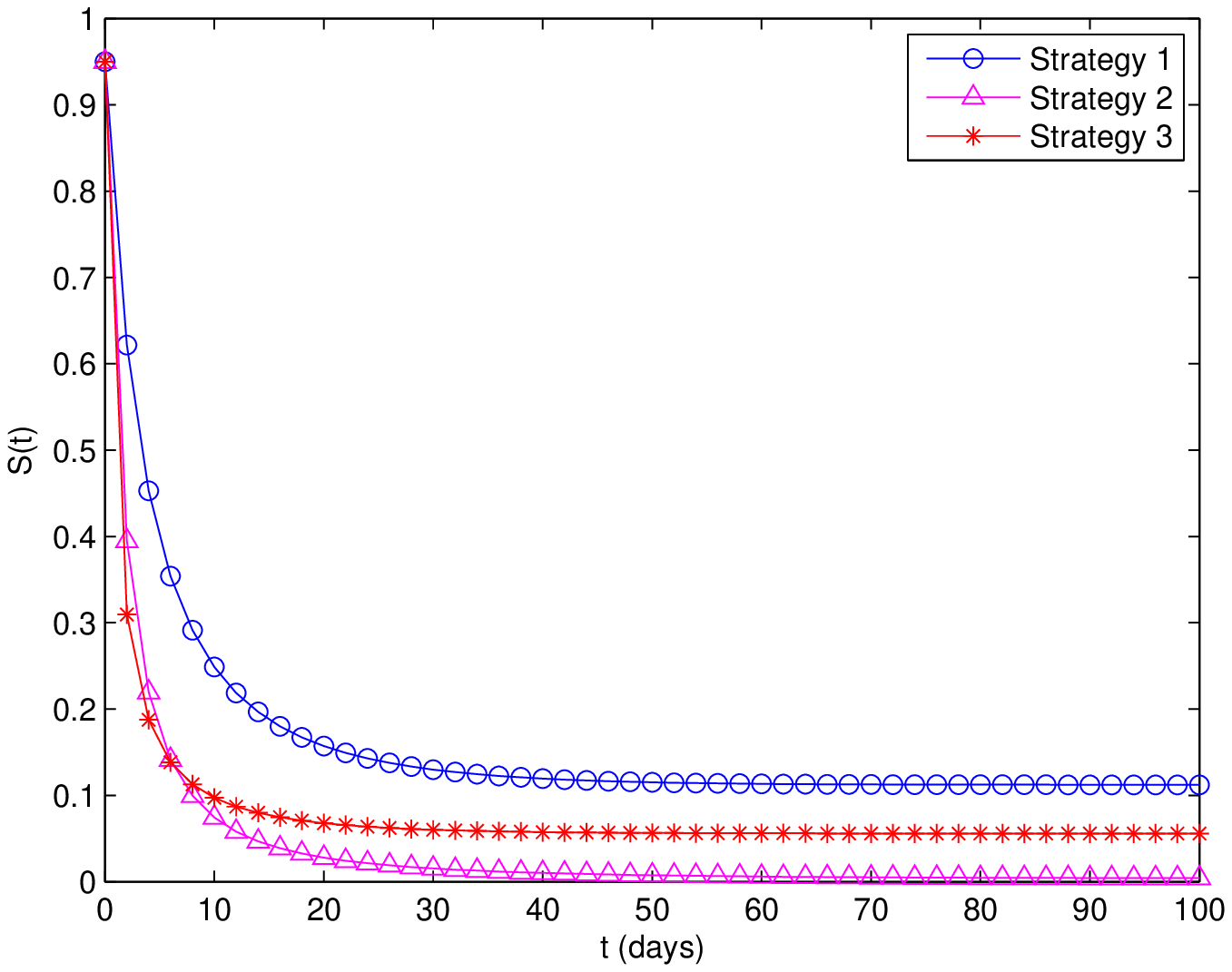}
\caption{Comparison between the curves of susceptible individuals $S(t)$
in case of Strategy~$1$, Strategy~$2$, and Strategy~$3$. \label{stratg123_S}}
\end{figure}

\begin{figure}
\centering
\includegraphics[scale=0.6]{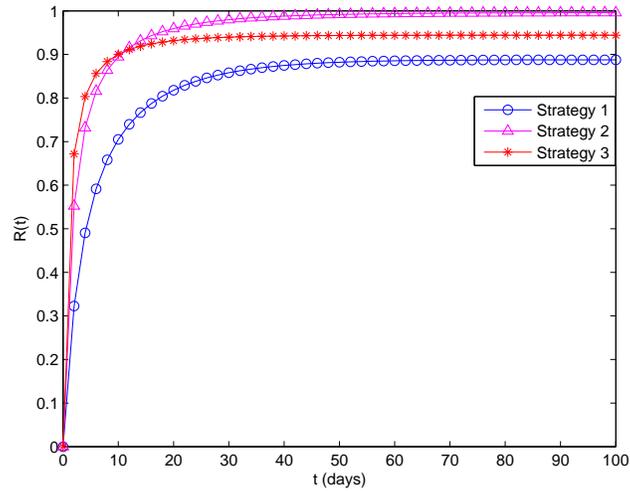}
\caption{Comparison between the curves of recovered individuals $R(t)$
in case of Strategy~$1$, Strategy~$2$, and Strategy~$3$. \label{stratg123_R}}
\end{figure}

\begin{figure}
\centering
\includegraphics[scale=0.6]{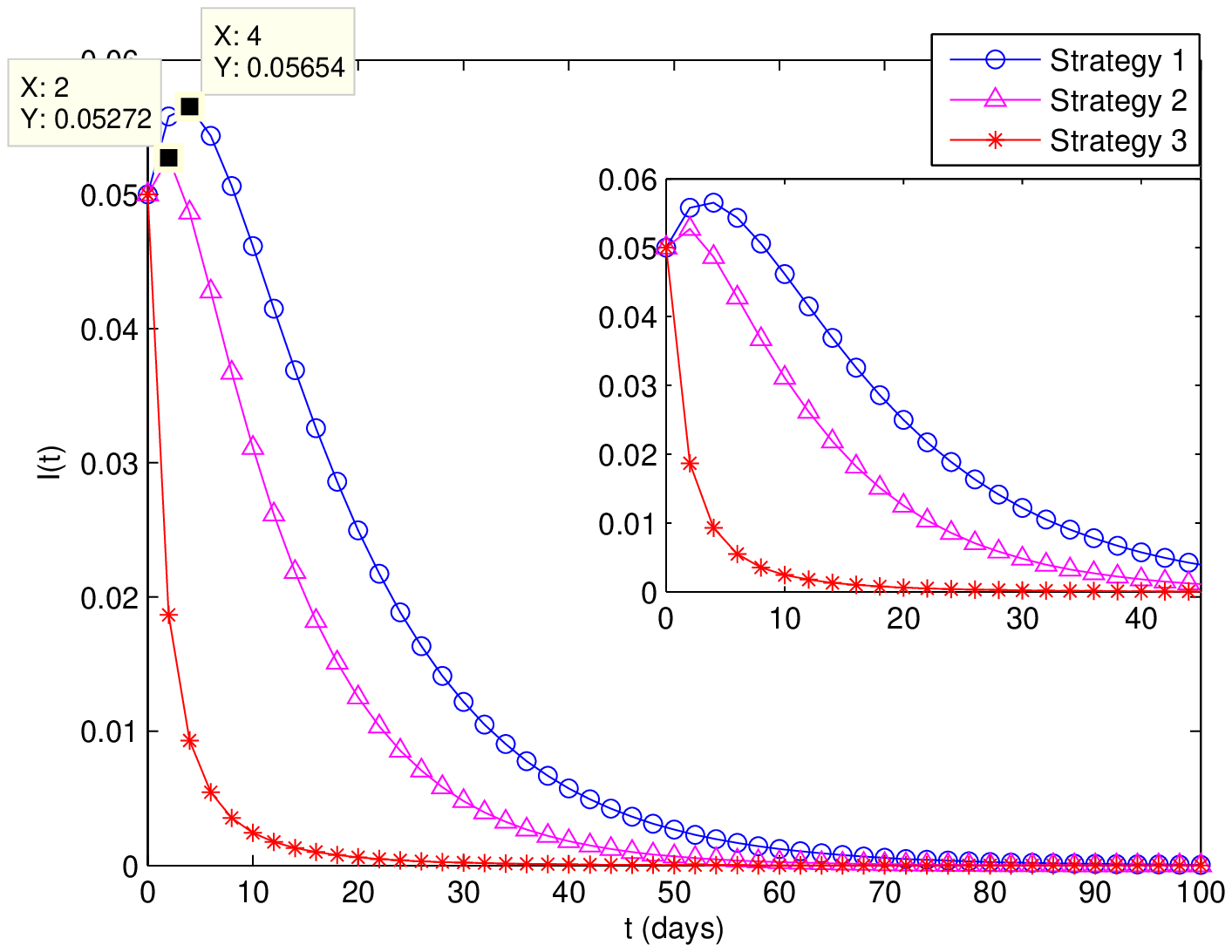}
\caption{Comparison between the curves of infected individuals $I(t)$
in case of Strategy~$1$, Strategy~$2$, and Strategy~$3$. \label{stratg123_I}}
\end{figure}

We compare between the three strategies
of Section~\ref{sec:3}, and we discuss the obtained results.
Figures~\ref{stratg123_S} and \ref{stratg123_R} represent, respectively, the
time-dependent curve of susceptible $S$ and recovered individuals $R$.
What is most important in the comparison
between the three curves, is the number of infected individuals that decreases
more rapidly in case of Strategy~$3$ (see  Figure~\ref{stratg123_I}).
Moreover, there isn't any peak in case of Strategy~$3$.
The results show the efficiency of an educational campaign
in controlling Ebola virus when we couple it with the treatment control
(as in Strategy~$3$). We conclude that one can improve results
by mixing vaccination and educational campaigns with treatment, which
has a great importance in poor countries that do not have the capacity
to defend themselves against the virus.
% -------------------------------
% -------------------------------
% -------------------------------

% ------------------------------------

\section{Conclusion}
\label{sec:4}

We improved the optimal control model discussed by Rachah and Torres
in \cite{MyID:321} (see also \cite{MyID:335}),
which provides a good description of the 2014
Ebola outbreak in West Africa. We discussed the integration of
an educational campaign about the virus into the population.
We have shown that the educational campaign has a great importance
with the treatment program, mainly in countries who don't have
the capacity to defend themselves against the virus. As future work,
we plan to include in our study other factors. For instance,
we intend to include in the mathematical model a quarantine procedure.

% -----------------------------------------

\subsection*{Acknowledgments}

This research was supported by the
Institut de Math\'{e}matiques de Toulouse, France (Rachah);
and by the Portuguese Foundation for Science and Technology (FCT),
within CIDMA project UID/MAT/04106/2013 (Torres).

% ------------------------------------------

% ------------------------------------

\end{document}